\documentclass[11 pt]{amsart}
\usepackage{hyperref}
\usepackage{ graphicx, mathrsfs}
\usepackage{amssymb,amsmath,amsthm}
\usepackage{mathrsfs,dsfont}

\theoremstyle{plain}
\newtheorem{theorem}{Theorem}

\newtheorem{proposition}{Proposition}
\newtheorem{lemma}[proposition]{Lemma}

\newtheorem{claim}{Claim}

\newtheorem*{remark}{Remark}
\begin{document}

\title{The second iterate for the Navier-Stokes equation}

\author{Pierre Germain}

\begin{abstract}
We consider the iterative resolution scheme for the Navier-Stokes equation, and focus on the second iterate, more precisely on the map from the initial data to the second iterate at a given time $t$.

We investigate boundedness properties of this bilinear operator. This new approach yields very interesting results: a new perspective on Koch-Tataru solutions; a first step towards weak strong uniqueness for Koch-Tataru solutions; and finally an instability result in $\dot{B}^{-1}_{\infty,q}$, for $q>2$.
\end{abstract}

\maketitle

\section{Introduction}

\subsection{The equation}

The Cauchy problem for the Navier-Stokes equation reads
$$
(NS) \;\;\;\;\left\{ \begin{array}{l} \partial_t u - \Delta u + {\mathbb{P}} \left( u \cdot \nabla u \right) = 0 \\
u_{|t=0} = u_0 \,\,, \end{array} \right.
$$
where $u$ is a function of $(t,x) \in [0,\infty) \times {\mathbb{R}}^d$ valued in ${\mathbb{R}}^d$. We denote $u\cdot \nabla = u^i \partial_i$ and $\mathbb{P}$ the Leray projector on divergence free vector fields, which is given on the Fourier side by $({\mathbb{P}} u)^\wedge (\xi) = P(\xi) \widehat{u}(\xi)$, with
$P_{ij}(\xi) = 1 - \frac{\xi_i \xi_j}{|\xi|^2}$.
 The initial data $u_0$ is taken to be divergence free. This feature is conserved by the flow of the above equation, namely $u$ remains divergence free at any later time.

\subsection{Iterative resolution of $(NS)$}

A natural way of solving the above system consists in setting up an iterative scheme: set $u_0 = 0$ and for any $n \geq 1$ let $u_n$ solve
\begin{equation}
\label{iteration1}
\left\{ \begin{array}{l} \partial_t u_n - \Delta u_n + {\mathbb{P}} \left( u_{n-1} \cdot \nabla u_{n-1} \right) = 0 \\
(u_{n})_{|t=0} = u_0 \,\,, \end{array} \right.
\end{equation}
If the sequence $(u_n)$ converges, the limit is, formally at least, a solution of $(NS)$ with initial data $u_0$.
Observe that $u_{n} - u_{n-1}$ is an $n$-linear operator from the data space to the solution space; let us denote it $F_n(u_0\,,\,\dots\,,\,u_0)$. Under appropriate convergence assumptions, one gets the following analytic expansion for the solution $u$
\begin{equation}
\label{expansionu}
u = \sum_{n=1}^{\infty} F_n (u_0\,,\,\dots\,,\,u_0)\,\,.
\end{equation}
The question is now to show convergence of the above iterative scheme; this is naturally done using functional spaces which are invariant by the scaling associated to $(NS)$
\begin{equation}
\label{scaling}
u_0(x) \longrightarrow \lambda u_0(\lambda x) \;\;\;\;\;\;\;\; u(x,t) \longrightarrow \lambda u(\lambda x\,,\,\lambda t) \,\,.
\end{equation}
This approach has been developped since the seminal paper of Fujita and Kato~\cite{FK} by (among others) Kato~\cite{FK}, Cannone~\cite{C}, and finally Koch and Tataru~\cite{KT}. These authors considered respectively data in the following scale-invariant spaces
$$
\dot{H}^{\frac{d}{2}-1} \hookrightarrow L^d \hookrightarrow \dot{B}^{-1+\frac{d}{p}}_{p,q} \hookrightarrow \nabla BMO\,\,\mbox{(with $d < p < \infty$ and $d \leq q \leq \infty$)} \,\,.
$$
Here we denoted $\dot{B}^{s}_{p,q}$ for the standard Besov space (see subsection~\ref{besov}), and $\nabla BMO$ for the space of derivatives of functions of $BMO$.
It is believed that $\nabla BMO$ is the optimal space in which the above scheme can be implemented.

The space $\dot{B}^{-1}_{\infty,\infty}$ is only a trifle larger than $\nabla BMO$, and contains any space of tempered distributions whose norm is invariant by~(\ref{scaling}). It has been showed by Montgomery-Smith~\cite{MS}, for an equation similar to $(NS)$, that the iterative scheme cannot be run in $\dot{B}^{-1}_{\infty,\infty}$. 

We prove in the present paper that the iterative scheme for $(NS)$ cannot be used for $u_0 \in \dot{B}^{-1}_{\infty,\infty}$; namely one cannot define properly the second iterate.

\subsection{$L^2$ stability of solutions}

Another approach to solving the above system is to use conservation of energy, and compactness arguments; one then obtains weak solutions in the energy class $L^\infty ([0,\infty),L^2) \cap L^2 ([0,\infty),\dot{H}^1)$, first constructed by Leray~\cite{L}. Is is not known in general whether these weak solutions are unique. 

In order to prove uniqueness of a given solution $w$ in the energy class, the weak-strong method essentially consists in proving $L^2$ stability for $(NS)$ around $w$. Due to the conservation of energy for $(NS)$, this essentially reduces to $L^2$ stability for the linearized equation around $w$. In other words, one needs to show that if $w$ solves
$$
\left\{ \begin{array}{l} \partial_t w - \Delta w + {\mathbb{P}} \left( w \cdot \nabla w \right) = 0 \\
w_{|t=0} = w_0 \,\,, \end{array} \right.
$$
the $L^2$ norm of $v$, solution of
\begin{equation}
\label{perteq}
\;\;\;\;\left\{ \begin{array}{l} \partial_t v - \Delta v + {\mathbb{P}} \left( v \cdot \nabla w + w \cdot \nabla v \right) = 0 \\
v_{|t=0} = v_0 \,\,, \end{array} \right.
\end{equation}
can be controlled by the $L^2$ norm of $v_0$ and some norm of $w_0$.

Such an $L^2$ stability result has been proved for $w_0$ belonging to one of the spaces in the following hierarchy
$$
\dot{B}^{-1+\frac{d}{p}}_{p,q} \hookrightarrow \dot{B}^{-1+\frac{d}{p}}_{L^{p,\infty},q} \hookrightarrow \dot{B}^{-1+\frac{d}{p}}_{M^{r,p},q} \hookrightarrow \dot{B}^{-1+\frac{d}{p}}_{M(\dot{H}^{d/p},L^2),q} \;\;\;\;\mbox{with}\;\;\frac{2}{q} + \frac{d}{p} = 1 \;\;,\;\;2 < r \leq p  \;\;\mbox{and}\;\; d < p < \infty \,\,. 
$$
(we refer to Subsection~\ref{besov} for a definition of the Besov spaces appearing above, and to the book of Lemari\'e-Rieusset~\cite{LR} for a definition of the Morrey spaces $M^{r,p}$ and of the multiplier spaces $M(\dot{H}^{d/q},L^2)$). 
The $L^2$ stability results for $w_0$ in the above spaces are due, respectively, to Prodi~\cite{P}, Dubois~\cite{D} and Lemari\'e-Rieusset~\cite{LR}.
The same result holds for $w_0$ small in one of the spaces in the following hierarchy
$$
L^d \hookrightarrow L^{d,\infty} \hookrightarrow M^{r,d} \hookrightarrow M(\dot{H}^{1},L^2) \;\;\;\;\mbox{with}\;\;2 < r \leq d \,\,,
$$
this is due respectively to Von Wahl~\cite{vonwahl}, Dubois~\cite{D} and Lemari\'e-Rieusset~\cite{LR}.

All the mentioned results were obtained by the following method: take the scalar product with $v$ of~(\ref{perteq}), and integrate in space and time. Optimal results using this method were obtained by the author~\cite{G} (it essentially consists of the spaces described above, with a broader range for $p$ and $q$). Thus a new method is needed if one wants to improve on the results just mentioned.

This is our aim here to follow a different approach. As in the previous subsection, we observe that $v$ can be expanded into a series of multilinear operators
\begin{equation}
\label{expansionv}
v = \sum_{n=0}^\infty G_{n}(v_0\,,\,w_0\,,\dots,\,w_0) \,\,,
\end{equation}
where $G_{n}(v_0\,,\,w_0\,,\dots,\,w_0)$ is linear in $v_0$ and $n$-linear in $w_0$. In this case though, a formula giving explicitly the $(G_n)$ would have to be more complicated than in the previous subsection.

We prove in the following that the second term (ie $n=1$) in the above expansion is bounded in $L^2$ if $v_0 \in L^2$, $w_0 \in \dot{B}^0_{\infty,\infty}$. We shall extend this result in a forthcoming article, and prove weak-strong uniqueness for Koch-Tataru solutions, that is solutions corresponding to data in $\nabla BMO$.

\section{Harmonic analysis background}

We adopt the unitary normalization for the Fourier transform, thus
$$
\widehat{f}(\xi) = \frac{1}{(2 \pi)^{d/2}} \int e^{-ix \xi} f(x) \, dx \;\;\;\mbox{and}\;\;\; f(x) = \frac{1}{(2 \pi)^{d/2}} \int e^{-ix \xi}\widehat{f}(\xi) \, d\xi \,\,.
$$

\subsection{Littlewood-Paley decomposition}

We define here the Littlewood-Paley operators, that will be of constant use in the following

To begin with, let us fix a dyadic partition of unity. In order to do so, pick a smooth function $\psi\;:\;\mathbb{R}^+ \rightarrow \mathbb{R}^+$ such that $\operatorname{Supp} \psi$ is a subset of $\left[ \frac{3}{4} \,;\, \frac{8}{3} \right]$ and
\begin{equation}
\label{psipsi}
\sum_{j \in \mathbb{Z}} \psi \left( \frac{\xi}{2^j} \right) = 1 \;\;\;\;\;\;\mbox{for $\xi > 0$} \,\,.
\end{equation}
Define the Fourier multipliers
\begin{equation*}
\begin{split}
&\Delta_N \overset{def}{=} \psi \left( \frac{|D|}{2^N} \right) \\
& \Delta_{\leq N} \overset{def}{=} \sum_{j \leq N} \psi \left( \frac{|D|}{2^j} \right) \\
& \Delta_{\geq N} \overset{def}{=} \sum_{j \geq N} \psi \left( \frac{|D|}{2^j} \right) \\
& \Delta_{M\leq\cdot\leq N} \overset{def}{=} \sum_{j = M}^N \psi \left( \frac{|D|}{2^j} \right) \,\,.
\end{split}
\end{equation*}
The two following identities hold (on any $L^p$, $1 \leq p < \infty$, and more generally on functional spaces which impose decay at infinity)
$$
\sum_j \Delta_j = \operatorname{Id} \;\;\;\;\mbox{and}\;\;\;\; \Delta_{\leq 0} + \sum_{j \geq 1} \Delta_j = 1\,\,.
$$

\subsection{Besov spaces}

\label{besov}

In this paragraph, we quickly define Besov spaces, and refer to the book of Lemari\'e-Rieusset~\cite{LR} for further information. Then we study the embedding properties of $\nabla BMO$ and $\dot{B}^{-1}_{\infty,p}$.

\bigskip

If $(s,q) \in \mathbb{R} \times [1,\infty] \times [1,\infty]$ and $E$ is a Banach space, the Besov space $\dot{B}^s_{E,q}$ is the space given by the norm
$$
\|f\|_{\dot{B}^{s}_{E,q}} \overset{def}{=} \left\| 2^{js} \left\| \Delta_j f \right\|_{E} \right\|_{\ell^q_j} = \left( \sum_{j\in \mathbb{Z}} 2^{jsq} \left\| \Delta_j f \right\|_{E}^{q} \right)^{1/q}
$$
(with the usual modification of the $\ell^q$ norms in case the index is infinite). This is a Banach space under appropriate conditions on $E$, $s$ and $q$. 

In case $E$ is the Lebesgue space $L^p$, with $1 \leq p \leq \infty$, we simpy denote $\dot{B}^{s}_{p,q}$ for $\dot{B}^{s}_{E,q}$. Thus
$$
\|f\|_{\dot{B}^{s}_{p,q}} \overset{def}{=} \left\| 2^{js} \left\| \Delta_j f \right\|_{L^p_x} \right\|_{\ell^q_j} = \left( \sum_{j\in \mathbb{Z}} 2^{jsq} \left\| \Delta_j f \right\|_{L^p}^{q} \right)^{1/q}
$$
In particular,
$$
\|f\|_{\dot{B}^{-1}_{\infty , q}} = \left( \sum_{j \in \mathbb{Z}} 2^{-jq} \|\Delta_j f\|_\infty^q \right)^{1/q} \;\;\;\;\mbox{for $1<q<\infty$, and} \;\;\;\;\|f\|_{\dot{B}^{-1}_{\infty , \infty}} = \sup_{j \in \mathbb{Z}} 2^{-j} \|\Delta_j f\|_\infty \,\,.
$$

\begin{lemma}
The following embeddings hold
\begin{equation*}
\begin{split}
& \dot{B}^{-1}_{\infty,q} \hookrightarrow \dot{B}^{-1}_{\infty,r} \;\;\;\;\;\;\mbox{if $q \leq r$} \\
& \dot{B}^{-1}_{\infty,2} \hookrightarrow \nabla BMO \hookrightarrow \dot{B}^{-1}_{\infty,\infty} \,\,.
\end{split}
\end{equation*}
\end{lemma}
\textsc{Proof:} The first embedding follows immediately from the definition of Besov spaces. The second embedding can be seen as an immediate consequence of the following equivalent characterization of the norms of $\dot{B}^{-1}_{\infty,2}$ and $BMO$ 
\begin{equation}
\label{eqBMO}
\begin{split}
& \|f\|_{\dot{B}^{-1}_{\infty,2}} \sim \int_0^\infty \|e^{t\Delta} f\|_{\infty}^2 \,dt \\
& \|f\|_{\nabla BMO} \sim \sup_{x,R} \frac{1}{R^d} \int_0^{R^2} \int_{B(x,R)} |(e^{t\Delta} f)(y) |^2 \,dy\,dt\,\,.
\end{split}
\end{equation}
As for the last embedding, it follows from the equivalent characterization of $\dot{B}^{-1}_{\infty,\infty}$
$$
\|f\|_{\dot{B}^{-1}_{\infty,\infty}} \sim \sup_{t>0} \sqrt{t} \|e^{t\Delta} f\|_{L^\infty} \lesssim \|f\|_{\nabla BMO}
$$
(the proof of the above equivalences can be found in the book of Lemari\'e-Rieusset~\cite{LR}). $\blacksquare$

\section{Main result}

Computing explicitly the two first terms of the expansion~(\ref{expansionu}) yields
$$
\widehat{u} (t,\xi) = e^{-t |\xi|^2} \widehat{u_0}(\xi) + P(\xi)e^{-t |\xi|^2} \int_0^t e^{s(|\xi|^2 - |\eta|^2 - |\xi -\eta|^2)} \left( \widehat{u_0}(\eta) \cdot (\xi - \eta) \widehat{u_0}(\xi - \eta) \right) \,d\eta\,ds + \dots \,\,.
$$
For the expansion~(\ref{expansionv}), one gets
\begin{equation*}
\begin{split}
\widehat{v} (t,\xi & ) =  e^{-t |\xi|^2} \widehat{v_0}(\xi) \\
& + P(\xi)e^{-t |\xi|^2} \int_0^t e^{s(|\xi|^2 - |\eta|^2 - |\xi -\eta|^2)}  \left( \widehat{v_0}(\eta) \cdot (\xi - \eta) \widehat{w_0}(\xi - \eta) + \widehat{w_0}(\eta) \cdot (\xi - \eta) \widehat{v_0}(\xi - \eta)\right) \,d\eta\,ds \\ & + \dots \,\,.
\end{split}
\end{equation*}
Our aim is to study the boundedness properties of the bilinear terms appearing in the above expressions. By scaling, it suffices to consider the case $t=1$. Thus we are interested in the boundedness properties of $T_1$, $T_2$ given by
\begin{equation*}
\begin{split}
&\left(T_1(f,g) \right)^\wedge (\xi) \overset{def}{=} P(\xi) e^{-|\xi|^2} \int_0^1 \int e^{s(|\xi|^2 - |\eta|^2 - |\xi -\eta|^2)} \left( \widehat{f}(\eta) \cdot (\xi - \eta) \widehat{g}(\xi - \eta) \right) \,d\eta\,ds \,\,. \\
&\left( T_2(f,g)\right)^\wedge (\xi) \overset{def}{=} \\
& \;\;\;\;\;\;\;\;P(\xi)e^{-|\xi|^2} \int_0^1 \int e^{s(|\xi|^2 - |\eta|^2 - |\xi -\eta|^2)}  \left( \widehat{f}(\eta) \cdot (\xi - \eta) \widehat{g}(\xi - \eta) + \widehat{g}(\eta) \cdot (\xi - \eta) \widehat{f}(\xi - \eta)\right) \,d\eta\,ds \,\,,
\end{split}
\end{equation*}
where $f$, $g$, as well as $T_i (f,g)$ are divergence free maps from $\mathbb{R}^d$ to $\mathbb{R}^d$.

\begin{theorem} 
\label{theoprincipal} (We denote in the statement of this theorem $E_\sigma$ for divergence free vector fields in the Banach space $E$. In the following, we shall drop the subscript $\sigma$ to make notations lighter.)

(i) The operator $T_2$ is bounded from $\left( \dot{B}^{-1}_{\infty,\infty} \right)_\sigma \times \left( L^2 \right)_\sigma$ to $\left( L^2 \right)_\sigma$.

\medskip

(ii) The operator $T_1$ is bounded from $\left( \nabla BMO \right)_\sigma \times \left( \nabla BMO \right)_\sigma$ to $\left( \nabla L^\infty \right)_\sigma$.

\medskip

(iii) The operator $T_1$ is not bounded from $\left( \dot{B}^{-1}_{\infty,q} \right)_\sigma \times \left( \dot{B}^{-1}_{\infty,q} \right)_\sigma$ to $\mathcal{S}'$ if $q>2$. This is the case even if one restricts $T_1$ to the diagonal $(f,f) \in \left( \dot{B}^{-1}_{\infty,q} \right)_\sigma \times \left( \dot{B}^{-1}_{\infty,q} \right)_\sigma$.

\end{theorem}

What does the above theorem mean for the PDE problems evoked in the introduction?

\begin{itemize}

\item Point $(i)$ of the theorem is a first step towards weak-strong uniqueness for the Koch-Tataru solutions; or in other words, towards proving that if the data $u_0$ belongs to $\nabla BMO \cap L^2$, then the solution to $(NS)$ is unique in the energy class.
We will address this question in a forthcoming article.

\item Point $(ii)$ of the theorem says that the second iterate of~(\ref{iteration1}) belongs to $L^\infty ([0,\infty),\nabla L^\infty )$ if $u_0 \in \nabla BMO$. This is not proving the existence of a solution for such data (theorem of Koch and Tataru), but we believe it does provide interesting insight: it corresponds to the frequency approach, whereas Koch and Tataru's proof is done in physical space. Notice that the regularity of the bilinear term, $L^\infty \nabla L^\infty$, is slightly better than that of the linear term, $L^\infty \nabla BMO$.

\item Point $(iii)$ proves that the iteration scheme cannot be used to build up solutions associated to data in $\dot{B}^{-1}_{\infty,q}$, if $q>2$. This is not quite an ill-posedness result, but it says that the map which associates to the data a solution of $(NS)$ cannot be of class $\mathcal{C}^2$ from $\dot{B}^{-1}_{\infty,q}$ to $\mathcal{S}'$.

\end{itemize}

We would like to mention here the article of Germain, Masmoudi and Shatah~\cite{GMS}. It somehow corresponds to the dispersive equation version of the approach followed here for $(NS)$, in particular in the use which is made of bilinear operators. The essential difference between these two settings is that (space-time) resonances, which are the key to understanding global existence for dispersive equations, are not so relevant for a dissipative equation.

Finally, we learned after completion of the present work from J. Bourgain and N. Pavlovi\v{c} that they have just finished the proof of an ill-posedness result for the Navier-Stokes equation in the spaces $\dot{B}^{-1}_{\infty,q}$, with $q>2$.

\section{Bilinear operators}

Before proceeding with the proof of Theorem~\ref{theoprincipal}, let us recall some classical facts about bilinear operators.

The pseudo-product operator $B$ associated to the symbol $m(\xi,\eta)$ is defined by
$$
\left( B_m (f,g) \right)^\wedge (\xi) = \int_{\mathbb{R}^d} m(\xi,\eta) \widehat{f}(\eta) \widehat{g}(\xi-\eta) \, d\eta\,\,.
$$
Under appropriate conditions, these bilinear operators enjoy the same bounds as the ones given by H\"older's inequality for the standard multiplication. This is the content of the celebrated theorem of Coifman and Meyer.

\begin{theorem}[Coifman-Meyer~\cite{coifmanmeyer}~\cite{kenigstein}~\cite{grafakostorres}] 
\label{theoCM}
If the symbol $m$ satisfies for sufficiently many derivatives
\begin{equation}
\label{estimatescm}
| \partial_\xi^\alpha \partial_\eta^\beta m(\xi,\eta)  | \lesssim \frac{1}{\left( |\xi| + |\eta| \right)^{|\alpha| + |\beta|}} \,\,,
\end{equation}
then the associated operator is bounded
$$
B_m : L^p \times L^q \longrightarrow L^r \;\;\;\;\;\mbox{with} \;\;\frac{1}{p} + \frac{1}{q} = \frac{1}{r}\;,\; 1<p,q\leq \infty\;,\; 1\leq r <\infty\,\,.
$$
\end{theorem}

Unfortunately, the above theorem misses the endpoint $(\infty,\infty,\infty)$. This shortcoming can be overcome by strengthening the conditions on $m$.

\begin{proposition}
\label{bilinlinf}
Suppose the symbol $m$ satisfies
$$
m\,,\,\left( \nabla_{\xi,\eta} \right)^{d + 1} m \;\in\;L^2_{\xi,\eta}\,\,.
$$
Then $T_m$ is bounded from $L^\infty \times L^\infty$ to $L^\infty$.
\end{proposition}

\textsc{Proof:} Notice first that, if $M$ denotes the inverse Fourier transform (in $\xi$ and $\eta$) of $m$, then
$$
B_m(f,g) = \int M(x-z\,,\,x+z-y) f(y) g(z) \,dy\,dz\,\,.
$$
Furthermore, under the assumptions of the proposition, $M \in L^1$. $\blacksquare$

\section{Proof of $(i)$ in Theorem~\ref{theoprincipal}} 

\subsection{Reduction of the problem}

\label{reduc}

First, it is possible to simplify a little the problem by observing that
\begin{itemize}
\item Since we are dealing with divergence free functions, it is possible to replace in the definition of $T_1$ and $T_2$ $\xi - \eta$ by $\xi$.
 \item One can forget the Leray projector at the beginning of the expression of $T_2$, since it reduces to Riesz projections, which are bounded on the spaces of interest for us.
\item The vectorial nature of the functions $f$ and $g$ in the definition of $T_2$ will not play any role. Thus we replace $f$ and $g$ by scalar functions.
\item A function in $\dot{B}^{-1}_{\infty,\infty}$ can be written $\sum_{i=1}^{d} \partial_i f^i$, with $f^i$ in $\dot{B}^{0}_{\infty,\infty}$.
\end{itemize}

The above considerations show that the boundedness of $T_2 : \dot{B}^{-1}_{\infty,\infty} \times L^2 \longrightarrow L^2$ is implied by the following

\begin{claim}
If $\ell^1$ and $\ell^2$ are linear functions, the operator $B_\mu$ with symbol
$$
\mu(\xi,\eta) = e^{-|\xi|^2}  \ell^1(\eta) \ell^2(\xi) \int_0^1 e^{s(|\xi|^2 - |\eta|^2 - |\xi -\eta|^2)}  \,ds
$$
is bounded from $\dot{B}^{0}_{\infty,\infty} \times L^2$ to $L^2$.
\end{claim}

\bigskip

The idea of the proof of this claim is to decompose the $(\xi,\eta)$ plane into three regions, by picking three smooth functions $\chi_1$, $\chi_2$ and $\chi_3$ of $(\xi,\eta)$ such that
\begin{equation*}
\begin{split} 
& \chi_1 + \chi_2 + \chi_3 = 1 \\
& \operatorname{Supp}(\chi_1) \subset \{ |\xi| + |\eta| \leq 2 \} \\
& \operatorname{Supp}(\chi_2) \subset \{ |\xi| + |\eta| \geq 1 \,,\, |\xi| \geq \frac{1}{6} |\eta| \} \\ 
& \operatorname{Supp}(\chi_3) \subset \{ |\xi| + |\eta| \geq 1 \,,\, |\xi| \leq \frac{1}{5} |\eta| \} \\
& \mbox{$\chi_2$ and $\chi_3$ homogeneous of degree 0 for $|\xi|+|\eta| \geq 3$} 
\end{split}
\end{equation*}

Now let us define $B_{\mu_1}$, $B_{\mu_2}$ and $B_{\mu_3}$ by their symbols
$$
\mu_i(\xi,\eta) \overset{def}{=}  \chi_i(\xi,\eta) \mu(\xi,\eta)\,\,.
$$
Then obviously $\mu_1 + \mu_2 + \mu_3 = \mu$.

\bigskip

In the following subsections, we prove the boundedness of $B_{\mu_1}$, $B_{\mu_2}$ and $B_{\mu_3}$ from $\dot{B}^{0}_{\infty,\infty} \times L^2$ to $L^2$. Since $B_\mu = B_{\mu_1} + B_{\mu_2} + B_{\mu_3}$, this shall prove the claim, hence part $(i)$ of Theorem~\ref{theoprincipal}.

\subsection{The region where $|\xi| + |\eta| \lesssim 1$: boundedness of $B_{\mu_1}:\dot{B}^{0}_{\infty,\infty} \times L^2 \rightarrow L^2$}

\label{moineau1}

It is clear that
$$
\Delta_{\leq 1} \ell^1 (D) \;: \; \dot{B}^{0}_{\infty,\infty} \longrightarrow L^\infty \;\;\;\mbox{is bounded}\,\,.
$$
We also observe that the symbol 
$$
\mu_1(\xi,\eta) (\ell^1(\eta))^{-1} = e^{-|\xi|^2} \chi_1(\xi,\eta) \ell^2(\xi) \int_0^1 e^{s(|\xi|^2 - |\eta|^2 - |\xi -\eta|^2)}  \,ds
$$ 
satisfies the conditions of Theorem~\ref{theoCM} (it actually even belongs to $\mathcal{C}^\infty_0$). Thus one can estimate
\begin{equation*}
\begin{split}
\| B_{\mu_1} (f,g) \|_2 & = \| B_{\mu_1} (\Delta_{\leq 1} f \,,\, g) \|_2 \\
& = \left\| B_{\mu_{1}(\xi,\eta)(\ell^1(\eta))^{-1}} \left( \Delta_{\leq 1} \ell^1 (D) f \,,\,g \right) \right\|_2 \\
& \lesssim \left\| \Delta_{\leq 1} \ell^1 (D) f \right\|_\infty \|g\|_2 \\
& \lesssim \|f\|_{\dot{B}^{0}_{\infty,\infty}} \|g\|_2 \,\,.
\end{split}
\end{equation*}

\subsection{The region where $|\xi| + |\eta| \gtrsim 1$ and $|\xi| \gtrsim |\eta|$: boundedness of $B_{\mu_2}:\dot{B}^{0}_{\infty,\infty} \times L^2 \rightarrow L^2$}

\label{moineau2}

On the one hand, 
$$
\Delta_{\geq -1} e^{\frac{1}{100}\Delta} \;: \; \dot{B}^{0}_{\infty,\infty} \longrightarrow L^\infty \;\;\;\mbox{is bounded}\,\,.
$$
On the other hand, the symbol 
$$
\mu_2(\xi,\eta) e^{\frac{1}{100}|\eta|^2} = e^{\frac{1}{100}|\eta|^2} e^{-|\xi|^2} \chi_2(\xi,\eta) \ell^2(\xi) \ell^1(\eta) \int_0^1 e^{s(|\xi|^2 - |\eta|^2 - |\xi -\eta|^2)}  \,ds
$$
satisfies the conditions of Theorem~\ref{theoCM}. Indeed, any derivative of this symbol decays like an inverse exponential of $|\xi|^2+|\eta|^2$. Let us see quickly how such a decay estimate can be obtained for $\mu_2(\xi,\eta) e^{\frac{1}{100}|\eta|^2}$, and it will become clear that the same holds for any derivative. Using first that $|\xi|^2 - |\eta|^2 - |\xi-\eta|^2 \leq \frac{1}{2} |\xi|^2$, and then that on the support of $\chi_2$, $|\xi| \geq \frac{1}{6}|\eta|$, we have
$$
|\mu_{2}(\xi,\eta) e^{\frac{1}{100}|\eta|^2}| \leq e^{\frac{1}{100}|\eta|^2} e^{-\frac{1}{2}|\xi|^2} \chi_2(\xi,\eta) |\ell^1(\eta) \ell^2(\xi)| \leq e^{-\frac{1}{10}|\xi|^2} \chi_2(\xi,\eta) |\ell^1(\eta) \ell^2(\xi)| \lesssim e^{-\frac{1}{1000}(|\xi|^2 + |\eta|^2)} \,\,.
$$
Thus we can estimate
\begin{equation*}
\begin{split}
\| B_{\mu_2} (f,g) \|_2 & = \| B_{\mu_2} (\Delta_{\geq -1} f \,,\, g) \|_2 \\
& = \left\| B_{\mu_2(\xi,\eta)e^{\frac{1}{100}|\eta|^2}} \left( \Delta_{\geq -1} e^{\frac{1}{100}\Delta} f \,,\,g \right) \right\|_2 \\
& \lesssim \left\| \Delta_{\geq -1} e^{\frac{1}{100}\Delta} f \right\|_\infty \|g\|_2 \\
& \lesssim \|f\|_{\dot{B}^{0}_{\infty,\infty}} \|g\|_2 \,\,.
\end{split}
\end{equation*}

\subsection{Further refinement in the region where $|\xi| + |\eta| \gtrsim 1$ and $|\xi| << |\eta|$}

\label{moineau3}

In this region, which essentially corresponds to the support of $\chi_3$, the idea is to integrate out in time and get
$$
\mu_3(\xi,\eta) = \chi_3(\xi,\eta) \frac{\ell^2(\xi) \ell^1(\eta)}{|\xi^2|-|\eta|^2-|\xi-\eta|^2} \left( e^{-|\xi-\eta|^2-|\eta|^2} - e^{-|\xi|^2} \right) \,\,.
$$
Thus we can decompose
$$
\mu_3(\xi,\eta) = \mu_3'(\xi,\eta) - \mu_3''(\xi,\eta) \;\;\;\;\mbox{where}\;\;\;\left\{ \begin{array}{l} \mu_3'(\xi,\eta) = \chi_3(\xi,\eta) \frac{\ell^2(\xi) \ell^1(\eta)}{|\xi^2|-|\eta|^2-|\xi-\eta|^2} e^{-|\xi-\eta|^2-|\eta|^2} \\ \mu_3''(\xi,\eta) = \chi_3(\xi,\eta) \frac{\ell^2(\xi) \ell^1(\eta)}{|\xi^2|-|\eta|^2-|\xi-\eta|^2} e^{-|\xi|^2} \,\,. \\ \end{array} \right.
$$

\subsection{Boundedness of $B_{\mu_3'}:\dot{B}^{0}_{\infty,\infty} \times L^2 \rightarrow L^2$}
Observing that 
$$
\mu_3'(\xi,\eta) e^{|\eta|^2} = \chi_3(\xi,\eta) \frac{\ell^2(\xi) \ell^1(\eta)}{|\xi^2|-|\eta|^2-|\xi-\eta|^2} e^{-|\xi-\eta|^2}
$$
satisfies the conditions of Theorem~\ref{theoprincipal}, we can easily estimate
\begin{equation*}
\begin{split}
\| B_{\mu_3'} (f,g) \|_2 & = \| B_{\mu_3'} (\Delta_{\geq -1}f\,,\,g) \|_2 \\
& = \| B_{{\mu_3'} e^{|\eta|^2}} ( \Delta_{\geq -1} e^{\Delta} f \,,\, g) \|_2 \\
& \lesssim \| \Delta_{\geq -1} e^{\Delta} f\|_\infty \|g\|_2 \\
& \lesssim \|f\|_{\dot{B}^{0}_{\infty,\infty}} \|g\|_2\,\,.
\end{split}
\end{equation*}

\subsection{Boundedness of $B_{\mu_3''}:\dot{B}^{0}_{\infty,\infty} \times L^2 \rightarrow L^2$}

\label{mioneau4}

Observe that the symbol
$$
\mu_3''(\xi,\eta) |\eta| = \chi_3(\xi,\eta) \ell^2(\xi) e^{-|\xi|^2} \frac{|\eta| \ell^1(\eta)}{|\xi^2|-|\eta|^2-|\xi-\eta|^2} \frac{1}{|\eta|}
$$
satisfies the conditions of Theorem~\ref{theoCM}. Using furthermore that 
$$
\Delta_{\geq -1} \frac{1}{|D|} \; : \; \dot{B}^0_{\infty,\infty} \rightarrow L^\infty
$$
is bounded, we can estimate
\begin{equation*}
\begin{split}
\| B_{\mu_3''} (f,g) \|_2 & = \| B_{\mu_3'} (\Delta_{\geq -1}f\,,\,g) \|_2 \\
& = \| B_{{\mu_3''} |\eta|} ( \Delta_{\geq -1} \frac{1}{|D|} f \,,\, g) \|_2 \\
& \lesssim \| \Delta_{\geq -1} \frac{1}{|D|} f\|_\infty \|g\|_2 \\
& \lesssim \|f\|_{\dot{B}^{0}_{\infty,\infty}} \|g\|_2\,\,.
\end{split}
\end{equation*}

\section{Proof of $(ii)$ in Theorem~\ref{theoprincipal}}

\label{sectioninf}

The proof of point $(ii)$ is similar to that of point $(i)$. For this reason, we only sketch it, but emphasize the modifications that need to be done.

\subsection{Reduction of the problem}

As in subsection~\ref{reduc}, we observe that the boundedness of $T_1 : \nabla BMO \times \nabla BMO \longrightarrow \nabla L^\infty$ is implied by the 

\begin{claim}
If $\ell^1$ and $\ell^2$ are linear functions, the operator $B_\nu$ with symbol
$$
\nu(\xi,\eta) = e^{-|\xi|^2} \ell^1(\eta) \ell^2(\xi-\eta) \int_0^1 e^{s(|\xi|^2 - |\eta|^2 - |\xi -\eta|^2)}  \,ds
$$
is bounded from $BMO \times BMO$ to $L^\infty$.
\end{claim}

As in subsection~\ref{reduc}, let us define $B_{\nu_1}$, $B_{\nu_2}$ and $B_{\nu_3}$ by their symbols
$$
\nu_i(\xi,\eta) \overset{def}{=}  \chi_i(\xi,\eta) \nu(\xi,\eta)\,\,.
$$
Then obviously $\nu_1 + \nu_2 + \nu_3 = \nu$.

In the following subsections, we prove the boundedness of $B_{\nu_1}$, $B_{\nu_2}$ and $B_{\nu_3}$ from $BMO \times BMO$ to $L^\infty$.

\subsection{The region where $|\xi| + |\eta| \lesssim 1$: boundedness of $B_{\nu_1}:BMO \times BMO \rightarrow L^\infty$}

We proceed as in subsection~\ref{moineau1}, using that
$$
\Delta_{\leq 1} \ell^1 (D) \,:\,BMO \longrightarrow L^\infty\;\;\;\mbox{is bounded}\,\,,
$$
and then Proposition~\ref{bilinlinf} instead of Theorem~\ref{theoCM}.

\subsection{The region where $|\xi| + |\eta| \gtrsim 1$ and $|\xi| \gtrsim |\eta|$: boundedness of $B_{\nu_2} : BMO \times BMO \rightarrow L^\infty$}

We proceed as in subsection~\ref{moineau2}, using that
$$
\Delta_{\geq -1} e^{\frac{1}{100}\Delta} \;: \; BMO \longrightarrow L^\infty \;\;\;\mbox{is bounded}\,\,.
$$
and then Proposition~\ref{bilinlinf} instead of Theorem~\ref{theoCM}.

\subsection{Further refinement in the region where $|\xi| + |\eta| \gtrsim 1$ and $|\xi| << |\eta|$}

We integrate out in time and get
$$
\nu_3(\xi,\eta) = \chi_3(\xi,\eta) \frac{\ell^1(\eta) \ell^2(\xi-\eta)}{|\xi|^2-|\eta|^2-|\xi-\eta|^2} \left( e^{-|\xi-\eta|^2-|\eta|^2} - e^{-|\xi|^2} \right) \,\,.
$$
Thus we can decompose
$$
\nu_3(\xi,\eta) = \nu_3'(\xi,\eta) - \nu_3''(\xi,\eta) \;\;\;\;
\mbox{where}\;\;\;\left\{ \begin{array}{l} \nu_3'(\xi,\eta) = \chi_3(\xi,\eta) \frac{\ell^1(\eta) \ell^2(\xi-\eta)}{|\xi^2|-|\eta|^2-|\xi-\eta|^2} e^{-|\xi-\eta|^2-|\eta|^2} \\ \nu_3''(\xi,\eta) = \chi_3(\xi,\eta) \frac{\ell^1(\eta) \ell^2(\xi-\eta)}{|\xi^2|-|\eta|^2-|\xi-\eta|^2} e^{-|\xi|^2} \,\,. \\ \end{array} \right.
$$

\subsection{Boundedness of $B_{\nu_3'}: BMO \times BMO \rightarrow L^\infty$}
Let us further decompose 
$$
N(\xi,\eta) \overset{def}{=} \chi_3(\xi,\eta) \frac{\ell^1(\eta) \ell^2(\xi-\eta)}{|\xi^2|-|\eta|^2-|\xi-\eta|^2}
$$
as
$$
N(\xi,\eta) = \sum_{j=1}^\infty N_j(\xi,\eta) \overset{def}{=} \sum_{j=1}^\infty \psi\left( \frac{|\eta|}{2^j} \right) N(\xi,\eta)
$$
(recall that $\psi$ is defined in~(\ref{psipsi})). Observe that $N_1$ satisfies the hypotheses of Proposition~\ref{bilinlinf}, hence $B_{N_1} \;:\;L^\infty \times L^\infty \longrightarrow L^\infty$ is bounded. By scaling, the $B_{N_j} \;:\;L^\infty \times L^\infty \longrightarrow L^\infty$ are uniformly bounded.

To prove boundedness of $B_{\nu_3'}$, we need the following lemma from Chemin~\cite{C}.

\begin{lemma}
There exists a constant $c$ such that
$$
\left\| e^{t\Delta} \Delta_j f \right\|_\infty \lesssim e^{-c 2^{2j}} \left\| \Delta_j f \right\|_\infty \,\,.
$$
\end{lemma}

We can conclude:
\begin{equation*}
\begin{split}
\| B_{\nu_3'} (f,g) \|_\infty & = \| B_{N} (e^{\Delta} f, e^{\Delta} g) \|_\infty \\
& \leq \sum_{j=0}^{\infty} \| B_{N_j} (\Delta_{j-1 \leq \cdot \leq j+1} e^{\Delta} f, \Delta_{j-1 \leq \cdot \leq j+1} e^{\Delta} g) \|_\infty \\
& \lesssim \sum_{j=0}^\infty \left\| \Delta_{j-1 \leq \cdot \leq j+1} e^{\Delta} f \right\|_\infty  \left\| \Delta_{j-1 \leq \cdot \leq j+1} e^{\Delta} g \right\|_\infty \\
& \lesssim \sum_{j=0}^\infty e^{-c 2^{2j} } \left\| \Delta_{j-1 \leq \cdot \leq j+1} f \right\|_\infty  \left\| \Delta_{j-1 \leq \cdot \leq j+1} g \right\|_\infty \\
& \lesssim \|f\|_{BMO} \|g\|_{BMO} \,\,.
\end{split}
\end{equation*}

\subsection{Boundedness of $B_{\nu_3''}: BMO \times BMO \rightarrow L^\infty$} 
An important feature of the symbol $\nu_3''$ is that it only allows frequency interactions of the type ``high - high gives low''. Examining a little more this symbol, it becomes clear that $B_{\nu_3''}$ can be written -up to an easily estimated operator-
$$
B_{\nu_3''}(f,g) = e^\Delta \sum_{j\geq -1\,,\,|j-k| \leq 1} B_{\alpha_j} (\Delta_j f,\Delta_j g) \,\,,
$$
where $\alpha$ belongs to $\mathcal{C}_0^\infty$ and $\alpha_j (\xi,\eta) = \alpha \left( \frac{\xi}{2^j} \,,\, \frac{\eta}{2^j} \right)$ (to make notations lighter, we consider in the following that $j =k$).

Expanding $\alpha$ in Fourier series, we see that for some constant $c$ depending on the support of $\alpha$,
$$
\alpha (\xi,\eta) = \sum_{m,n \in \mathbb{Z}^d} \lambda_{m,n} e^{i c \left(m \eta + n (\xi - \eta) \right) }\,\,.
$$
This implies that 
$$
B_{\nu_3''}(f,g)(x) = e^\Delta \sum_{j \geq -1} \sum_{m,n \in \mathbb{Z}^d} \lambda_{m,n} \Delta_j f(x + c 2^{-j} m) \Delta_j g(x + c 2^{-j} n) \,\,.
$$
We can assume that the norms of $f$ and $g$ in $BMO$ are comparable. Taking advantage of the strong decay of the kernel associated to $e^\Delta$, we can estimate
\begin{equation*}
\begin{split}
&\left| e^\Delta \sum_{j \geq -1} \sum_{m,n \in \mathbb{Z}^d} \lambda_{m,n} \Delta_j f(x + c 2^{-j} m) \Delta_j g(x + c 2^{-j} n) \right| \\
 &\;\;\;\;\;\;\;\;\;\;\;\;\;\;\;\;\;\;\lesssim sup_{x \in \mathbb{R}^d} \sum_{m,n \in \mathbb{Z}^d} \lambda_{m,n}\left\| \sum_{j \geq -1} \Delta_j f(x + c 2^{-j} m) \Delta_j g(x + c 2^{-j} n)  \right\|_{L^1(B(x,1))} \\
& \;\;\;\;\;\;\;\;\;\;\;\;\;\;\;\;\;\; \lesssim sup_{x \in \mathbb{R}^d} \sum_{m,n \in \mathbb{Z}^d} \lambda_{m,n}\int_{B(x,1)} \left[ \sum_{j \geq -1} \left| \Delta_j f(x + c 2^{-j} m) \right|^2 + \left| \Delta_j g(x + c 2^{-j} n) \right|^2 \right]\,dx\,\,.
\end{split}
\end{equation*}
Using the following characterization of the norm of $BMO$, which is essentially a rephrasing of~(\ref{eqBMO}),
$$
\| f \|_{BMO} = \sup_{J \in \mathbb{Z}\,,\, x \in \mathbb{R}^d} \frac{1}{2^{Jd}} \int_{B(x,2^{J})} \sum_{j \geq -J} |\Delta_j f|^2 \,dx\,\,.
$$
yields now
\begin{equation*}
\begin{split}
\left| e^\Delta \sum_{j \geq -1} \sum_{m,n \in \mathbb{Z}^d} \lambda_{m,n} \Delta_j f(x + c 2^{-j} m) \Delta_j g(x + c 2^{-j} n) \right| &\lesssim \sum_{m,n \in \mathbb{Z}^d} \lambda_{m,n} (|m|+1)^d (|n|+1)^d \|f\|_{BMO} \|g\|_{BMO} \\
& \lesssim \|f\|_{BMO} \|g\|_{BMO}
\end{split}
\end{equation*}
where we used in the last inequality the rapid decay of the $(\lambda_{m,n})$.

\section{Proof of $(iii)$ in Theorem~\ref{theoprincipal}}

If $q>2$, we want to build up a counterexample to boundedness of
$$
T_2 : \dot{B}^{-1}_{\infty,q} \times \dot{B}^{-1}_{\infty,q} \longrightarrow \mathcal{S}'\,\,.
$$

\subsection{Idea behind the counterexample} 

Examining the analysis performed in the preceding section, it appears that if one excludes the region $|\xi| + |\eta| \gtrsim 1$, $|\xi| << |\eta|$, the operator $T_2 : \dot{B}^{-1}_{\infty,\infty} \times \dot{B}^{-1}_{\infty,\infty} \rightarrow \nabla L^\infty$ is bounded.

Thus our example should generate a ``high - high gives low'' frequency interaction, which becomes unbounded.

This insight is actually the only use that we shall make of the preceding analysis: the counterexample will be otherwise self-contained.

\subsection{The counterexample}

For simplicity, we set $d=3$ and pick $e_1 = \left( \begin{array}{l} 1 \\ 0 \\ 0 \end{array} \right)$, $e_2 = \left( \begin{array}{l} 0 \\ 1 \\ 0 \end{array} \right)$, $e_3 = \left( \begin{array}{l} 0 \\ 0 \\ 1 \end{array} \right)$ an orthonormal basis of $\mathbb{R}^3$. We shall denote $\times$ the vector product.

Next, let $\phi$ be a smooth, even, non-negative (real-valued) function on $\mathbb{R}^d$, such that $\phi = 1$ on $B(0,2)$ and $\phi = 0$ outside of $B(0,3)$. Also let $(\alpha_k)$ be a sequence in $\ell^q \setminus \ell^2$.
Define $f^N$ by its Fourier transform
\begin{equation*}
\begin{split}
\widehat{f^N}(\xi) & = \sum_{k=10}^N \widehat{f_{k,+}}(\xi) - \sum_{k=10}^N \widehat{f_{k,-}}(\xi) \\
& \overset{def}{=} \sum_{k=10}^N 2^k \alpha_k \phi (\xi - 2^k e_1) \left( \frac{\xi}{|\xi|} \times e_2 \right) - \sum_{k = 10}^N 2^k \alpha_k \phi (\xi + 2^k e_1) \left( \frac{\xi}{|\xi|} \times e_2 \right) \,\,.
\end{split}
\end{equation*}
It is clear that $f^N$ is real-valued, divergence-free, and uniformly bounded (with respect to $N$) in $\dot{B}^{-1}_{\infty,q}$. 

\begin{remark}
Let us pause for a moment and make two observations
\begin{itemize}
 \item First, notice that the above sequence is very similar to the one used by Montgomery-Smith~\cite{MS} to prove the result mentioned in the introduction, namely that for a Navier-Stokes like equation the iterative resolution method does not work for data in $\dot{B}^{-1}_{\infty,\infty}$. This is also very similar to the example used by Stein~\cite{S} to prove that symbols in $S^0_{1,1}$ are not in general associated to operators which are bounded on $L^2$. Thus, as Montgomery-Smith puts it, it might be that the non-boundedness result which we are about to prove ``says more about the nature of the $\dot{B}^{-1}_{\infty,\infty}$ space than about the Navier-Stokes equation itself''.
\item Second, we believe it is very instructive to relate the instability result for the data $f^N$ - that we will momentarily prove - to a result proved by Chemin and Gallagher~\cite{CG}. These authors build up data which are large in $\dot{B}^{-1}_{\infty,\infty}$ but still yield global solutions of the Navier-Stokes equation. These data have the following peculiarity: the different scales are physically separated, in other words the oscillations at different scales occur at different places; this is ensured by a fractal like transformation. This is to be contrasted with the $(f^N)$ for which oscillations at \textit{all} scales occur at the same location. 
\end{itemize}

\end{remark}

From now on, we fix 
$$
\xi_0 \overset{def}{=} \left( \begin{array}{l} 0 \\ \frac{1}{2} \\ \frac{1}{2} \end{array} \right)\;\;\;\;\mbox{thus}\;\;\;\; P(\xi_0) = 
\left( \begin{array}{ccc} 1 & 0 & 0 \\ 0 & \frac{1}{2} & -\frac{1}{2} \\ 0 & -\frac{1}{2} & \frac{1}{2} \end{array} \right) \,\,.
$$
An important and elementary observation is that the only possible interaction of $f^N$ with itself yielding this frequency $\xi_0$ corresponds to $f_{k,\pm}$ interacting with $f_{k,\mp}$.

This observation, along with performing the time integral in the definition of $T_1$, yields

\begin{equation*}
\begin{split}
T_1(f^N,f^N)^\wedge(\xi_0) & = P(\xi_0) e^{-|\xi_0|^2} \int \frac{1-e^{|\xi_0|^2 - |\eta|^2 - |\xi_0-\eta|^2}}{|\xi_0|^2 - |\eta|^2 - |\xi_0-\eta|^2} \xi_0 \cdot \left( \frac{\eta}{|\eta|} \times e_2 \right) \left( \frac{\xi_0-\eta}{|\xi_0-\eta|} \times e_2 \right)  \\
& \;\;\;\;\;\;\;\;\;\;\;\;\;\;\;\;\;\;\;\; \sum_{k=10}^N \alpha_k^2 2^{2k} \phi(\eta - 2^k e_1) \phi(\xi_0 - \eta + 2^k e_1) \,d\eta\\
+ & P(\xi_0) e^{-|\xi_0|^2} \int \frac{1-e^{|\xi_0|^2 - |\eta|^2 - |\xi_0-\eta|^2}}{|\xi_0|^2 - |\eta|^2 - |\xi_0-\eta|^2} \xi_0 \cdot \left( \frac{\eta}{|\eta|} \times e_2 \right) \left( \frac{\xi_0-\eta}{|\xi_0-\eta|} \times e_2 \right)  \\
& \;\;\;\;\;\;\;\;\;\;\;\;\;\;\;\;\;\;\;\; \sum_{k=10}^N \alpha_k^2 2^{2k} \phi(\eta + 2^k e_1) \phi(\xi_0 - \eta - 2^k e_1) \,d\eta
\end{split}
\end{equation*}
or, reorganizing things a little,
\begin{equation*}
\begin{split}
T_1(f^N,f^N)^\wedge  (\xi_0)  = & \int e^{-|\xi_0|^2}  \phi(\eta) \phi(\xi_0 - \eta) \\
& \left[ \sum_{k=10}^N \alpha^2_k \frac{2^{2k}(1 - e^{|\xi_0|^2 - |\eta+2^k e_1|^2 - |\xi_0-\eta - 2^k e_1|^2})}{ |\xi_0|^2 - |\eta+2^k e_1|^2 - |\xi_0-\eta - 2^k e_1|^2 } \right.\\
& \;\;\;\;\;\;\;\;\;\;\;\;\;\;\;\;\;\;\;\; P(\xi_0) \left[ \xi_0 \cdot  \left( \frac{\eta + 2^k e_1}{|\eta + 2^k e_1|} \times e_2 \right)\right] \left( \frac{\xi_0-\eta-2^k e_1}{|\xi_0-\eta-2^k e_1|} \times e_2 \right)   \\
& +  \sum_{k=10}^N \alpha^2_k \frac{2^{2k}(1 - e^{|\xi_0|^2 - |\eta+2^k e_1|^2 - |\xi_0-\eta - 2^k e_1|^2})}{ |\xi_0|^2 - |\eta - 2^k e_1|^2 - |\xi_0 - \eta + 2^k e_1|^2 } \\
& \;\;\;\;\;\;\;\;\;\;\;\;\;\;\;\;\;\;\;\;  P(\xi_0) \left.\left[  \xi_0 \cdot
\left( \frac{\eta - 2^k e_1}{|\eta - 2^k e_1|} \times e_2 \right) \right]\left( \frac{\xi_0-\eta + 2^k e_1}{|\xi_0-\eta  + 2^k e_1|} \times e_2 \right)  \right] \, d\eta 
\end{split}
\end{equation*}

It is easily seen that if $\eta \in \operatorname{Supp} \phi$, and $k \geq 10$
\begin{equation*}
\begin{split}
& \frac{2^{2k}(1 - e^{|\xi_0|^2 - |\eta+2^k e_1|^2 - |\xi_0-\eta - 2^k e_1|^2})}{ |\xi_0|^2 - |\eta+2^k e_1|^2 - |\xi_0-\eta - 2^k e_1|^2 } \left( \frac{\eta + 2^k e_1}{|\eta + 2^k e_1|} \times e_2 \right) \left( \frac{\xi_0-\eta-2^k e_1}{|\xi_0-\eta-2^k e_1|} \times e_2 \right) \sim - \frac{1}{2} \left( \begin{array}{l} 0 \\ 0 \\ 1 \end{array} \right) \otimes \left( \begin{array}{l} 0 \\ 0 \\ 1 \end{array} \right)\\
& \frac{2^{2k}(1 - e^{|\xi_0|^2 - |\eta+2^k e_1|^2 - |\xi_0-\eta - 2^k e_1|^2})}{ |\xi_0|^2 - |\eta - 2^k e_1|^2 - |\xi_0-\eta + 2^k e_1|^2 } \left( \frac{\eta - 2^k e_1}{|\eta - 2^k e_1|} \times e_2 \right) \left( \frac{\xi_0-\eta + 2^k e_1}{|\xi_0-\eta +2^k e_1|} \times e_2 \right) \sim - \frac{1}{2} \left( \begin{array}{l} 0 \\ 0 \\ 1 \end{array} \right) \otimes \left( \begin{array}{l} 0 \\ 0 \\ 1 \end{array} \right) \,\,.
\end{split}
\end{equation*}
We conclude that if $\eta \in \operatorname{Supp} \phi$, and $k \geq 10$
\begin{equation*}
\begin{split}
& \frac{2^{2k}(1 - e^{|\xi_0|^2 - |\eta+2^k e_1|^2 - |\xi_0-\eta - 2^k e_1|^2})}{ |\xi_0|^2 - |\eta+2^k e_1|^2 - |\xi_0-\eta - 2^k e_1|^2 } P(\xi_0) \left[ \xi_0 \cdot  \left( \frac{\eta + 2^k e_1}{|\eta + 2^k e_1|} \times e_2 \right)\right] \left( \frac{\xi_0-\eta-2^k e_1}{|\xi_0-\eta-2^k e_1|} \times e_2 \right)
 \sim \left( \begin{array}{l} 0 \\ \frac{1}{8} \\ -\frac{1}{8} \end{array} \right) \\
& \frac{2^{2k}(1 - e^{|\xi_0|^2 - |\eta+2^k e_1|^2 - |\xi_0-\eta - 2^k e_1|^2})}{ |\xi_0|^2 - |\eta - 2^k e_1|^2 - |\xi_0-\eta + 2^k e_1|^2 } P(\xi_0) \left[  \xi_0 \cdot
\left( \frac{\eta - 2^k e_1}{|\eta - 2^k e_1|} \times e_2 \right) \right]\left( \frac{\xi_0-\eta + 2^k e_1}{|\xi_0-\eta + 2^k e_1|} \times e_2 \right)\sim \left( \begin{array}{l} 0 \\ \frac{1}{8} \\ -\frac{1}{8} \end{array} \right) \,\,.
\end{split}
\end{equation*}
Integrating over $\eta$ and taking advantage of the positivity of $\phi$, we see that there exists a constant $C \neq 0$ such that
$$
T_1(f^N,f^N)^\wedge  (\xi_0) \sim C \sum_{k=1}^N \alpha_k^2 \left( \begin{array}{l} 0 \\ 1 \\ -1 \end{array} \right) \,\,,
$$
and in particular
$$
\left( T_1(f^N,f^N)^\wedge \right)^3 (\xi_0) \gtrsim \sum_{k=1}^N \alpha_k^2\,\,.
$$
We can run the same argument in a neighbourhood of $\xi_0 = \left( \begin{array}{l} 0 \\ \frac{1}{2} \\ \frac{1}{2} \end{array} \right)$, and obtain that, uniformly in $\zeta \in B\left( \xi_0 \,,\,\epsilon \right)$, for $\epsilon$ small enough,
$$
\left| \left( T_1(f^N,f^N)^\wedge \right)^3 (\zeta) \right| \gtrsim \sum_{k=10}^N \alpha_k^2\,\,.
$$
The series in the right hand side diverges. Thus, in spite of the boundedness of $f^N$ in $\dot{B}^{-1}_{\infty,q}$, the Fourier transform of $T_1(f^N,f^N)$ is, on $B\left( \xi_0 \,,\,\epsilon \right)$, larger than a diverging sequence. This means that $T_1(f^N,f^N)$ is not bounded in $\mathcal{S}'$. $\blacksquare$

\end{document}